\input amstex
\documentstyle{amsppt}
%
\catcode`@=11
\redefine\output@{%
  \def\break{\penalty-\@M}\let\par\endgraf
  \ifodd\pageno\global\hoffset=105pt\else\global\hoffset=8pt\fi  
  \shipout\vbox{%
    \ifplain@
      \let\makeheadline\relax \let\makefootline\relax
    \else
      \iffirstpage@ \global\firstpage@false
        \let\rightheadline\frheadline
        \let\leftheadline\flheadline
      \else
        \ifrunheads@ 
        \else \let\makeheadline\relax
        \fi
      \fi
    \fi
    \makeheadline \pagebody \makefootline}%
  \advancepageno \ifnum\outputpenalty>-\@MM\else\dosupereject\fi
}
\def\Beta{\mathchar"0\hexnumber@\rmfam 42}
\catcode`\@=\active
\nopagenumbers
\chardef\textvolna='176

\chardef\bigalpha='013
\def\negskp{\hskip -2pt}

\def\sign{\operatorname{sign}}

\def\compos{\,\raise 1pt\hbox{$\sssize\circ$} \,}

\def\inf{\operatornamewithlimits{inf}}

\def\blue#1{#1}

\def\registered{\ooalign {\hfil \raise .07ex\hbox {$\sssize\text{R}$}\hfil\crcr\Orb }}
\catcode`#=11\def\diez{#}\catcode`#=6
\catcode`&=11\catcode`&=4
\catcode`_=11\catcode`_=8
\catcode`~=11\catcode`~=\active
\def\mycite#1{\cite{\blue{#1}}\immediate\special{ps:
     ShrHPSdict begin /ShrBORDERthickness 0 def}}

\def\mytag#1{%
    \tag#1}
\def\mythetag#1{\thetag{\blue{#1}}\immediate\special{ps:
     ShrHPSdict begin /ShrBORDERthickness 0 def}}
\def\myrefno#1{\no#1}
\def\myhref#1#2{\blue{#2}\immediate\special{ps:
     ShrHPSdict begin /ShrBORDERthickness 0 def}}
\def\myEarXivlink{\myhref{http://arXiv.org}{http:/\negskp/arXiv.org}}

\def\mytheorem#1{\csname proclaim\endcsname{Theorem #1}}
\def\mytheoremwithtitle#1#2{\csname proclaim\endcsname{Theorem #1#2}}

\def\mylemma#1{\csname proclaim\endcsname{Lemma #1}}
\def\mylemmawithtitle#1#2{\csname proclaim\endcsname{Lemma #1#2}}

\def\mycorollary#1{\csname proclaim\endcsname{Corollary #1}}

\def\myconjecture#1{\csname proclaim\endcsname{Conjecture #1}}
\def\myconjecturewithtitle#1#2{\csname proclaim\endcsname{Conjecture #1#2}}

\def\myproblem#1{\csname proclaim\endcsname{Problem #1}}
\def\myproblemwithtitle#1#2{\csname proclaim\endcsname{Problem #1#2}}


\pagewidth{360pt}
\pageheight{606pt}
\topmatter
\title
On root mean square approximation 
by exponential functions.
\endtitle
\rightheadtext{On root mean square approximation \dots}
\author
Ruslan Sharipov
\endauthor
\address Bashkir State University, 32 Zaki Validi street, 450074 Ufa, Russia
\endaddress
\email\myhref{mailto:r-sharipov\@mail.ru}{r-sharipov\@mail.ru}
\endemail
\abstract
     The problem of root mean square approximation of a square integrable 
function by finite linear combinations of exponential functions is considered. 
It is subdivided into linear and nonlinear parts. The linear approximation 
problem is solved. Then the nonlinear problem is studied in some particular 
example. 
\endabstract
\endtopmatter
\TagsOnRight
\document

\head
1. Introduction.
\endhead
     Exponential functions of the form $e^{\kern 0.5pt\lambda\,x}$ arise as
solutions of linear differential equations. In physics they describe
oscillatory and damped oscillatory processes. Assume that $f(x)$ is a 
complex-valued square integrable function with the argument $x$ belonging to 
the interval $[\kern 1pt-\pi,\,+\pi]$ of the real line $\Bbb R$:
$$
\hskip -2em
f(x)\in L^2([\kern 1pt-\pi,\,+\pi]).
\mytag{1.1}
$$ 
Let's consider a finite sequence of exponential functions 
$$
\hskip -2em
e^{\kern 0.5pt\lambda_1x},\,\ldots,\,e^{\kern 0.5pt\lambda_nx},
\mytag{1.2}
$$
where $\lambda_1,\,\ldots,\,\lambda_n$ are distinct complex numbers, i\.\,e\.
$\lambda_{\kern 1pt i}\neq\lambda_j$. Using the functions \mythetag{1.2}, we
compose a linear combination with complex coefficients:
$$
\hskip -2em
\phi(x)=\sum^n_{i=1}a_i\,e^{\kern 0.5pt\lambda_ix}. 
\mytag{1.3}
$$
We say that the function \mythetag{1.3} approximates the square integrable 
function \mythetag{1.1} if the $L^2$-norm of their difference is sufficiently
small: 
$$
\hskip -2em
\Vert\kern 1ptf-\phi\kern 1pt\Vert=\sqrt{\frac{1}{2\,\pi}
\int\limits^{\,+\pi}_{\!-\pi}\!\!|\kern 1pt f(x)-\phi(x)\kern 1pt|^2
\,d\kern 0.5pt x}.
\mytag{1.4}
$$
The quantity \mythetag{1.4} is also known as the root mean square deflection 
of $\phi$ from $f$. The problem of minimizing this deflection is called the 
root mean square approximation problem. It can be attributed to the class 
of variational problems.\par 
     The problem of minimizing the root mean square deflection \mythetag{1.4}
is subdivided into linear and nonlinear parts. The linear problem consists in
finding optimal coefficients $a_1,\,\ldots,\,a_n$ in \mythetag{1.3}, provided
$\lambda_1,\,\ldots,\,\lambda_n$ are fixed. This problem is similar to those
studied by A. F. Leontiev and his school (see \mycite{1}). In the case of a
finite set of exponential functions \mythetag{1.2} it is solved completely in
an explicit form.\par
     The nonlinear approximation problem consists in minimizing the solution 
of the linear problem by varying $\lambda_1,\,\ldots,\,\lambda_n$ and choosing 
optimal values for them. It arises from the applied problem of numerical 
separation of a noised signal presumably being a mixture of oscillatory and 
damped oscillatory signals. In this form the problem was suggested by 
A.~S.~Vishnevskiy, president of PhysTech Co\., the weighing technologies 
company. In the present paper the nonlinear problem is 
studied in the example of the very simple function $f(x)=\sign(x)$.\par
\head
2. Solution of the linear approximation problem. 
\endhead
     From \mythetag{1.4} one can easily derive the following formula for the
the root mean square deflection of the function \mythetag{1.3} from $f$:
$$
\hskip -2em
\gathered
\Vert f-\phi\Vert^2=\Vert f\Vert^2-\sum^n_{j=1}a^{\kern 0.5pt j}\,\bigl<f\bigl|
e^{\kern 0.5pt\lambda_jx}\bigr>\,-\\
\kern 7em-\sum^n_{i=1}\overline{a^i}\,\bigl<e^{\kern 0.5pt\lambda_ix}\bigl|f\bigr>
+\sum^n_{i=1}\sum^n_{j=1}g_{ij}\,\overline{a^i}\,a^j.
\endgathered
\mytag{2.1}
$$
By means of angular brackets in \mythetag{2.1} we denote the $L^2$-scalar 
product 
$$
\hskip -2em
\bigl<a\bigl|b\bigr>=\frac{1}{2\,\pi}\int\limits^{\,+\pi}_{\!-\pi}
\overline{a(x)}\,b(x)\,d\kern 0.5pt x
\mytag{2.2}
$$
Overlined variables and functions in \mythetag{2.1} and \mythetag{2.2} mean
complex conjugates. Through $g_{ij}$ in \mythetag{2.1} we denote the components
of the Gram matrix $G$:
$$
\hskip -2em
g_{ij}=\bigl<e^{\kern 0.5pt\lambda_ix}\bigl|\,e^{\kern 0.5pt\lambda_jx}\bigr>.
\mytag{2.3}
$$
\par
     From \mythetag{2.2} one can easily derive the following property of
the $L^2$-scalar product 
$$
\hskip -2em
\overline{\bigl<a\bigl|\vphantom{\vrule height 10pt depth 0pt}b\bigr>}
=\bigl<b\bigl|a\bigr>.
\mytag{2.4}
$$
The property \mythetag{2.4} implies the following relationships:
$$
\xalignat 2
&\hskip -2em
g_{ij}=\overline{g_{j\kern 0.5pt i}\vphantom{\vrule height 6pt depth 0pt}},
&&\bigl<f\bigl|\kern 1pt e^{\kern 0.5pt\lambda_ix}\bigr>=
\overline{\bigl<\kern 0.5pt e^{\kern 0.5pt\lambda_ix}\bigl|f\bigr>
\vphantom{\vrule height 10pt depth 0pt}}.
\mytag{2.5}
\endxalignat
$$
\par 
     If we denote $F=\Vert f-\phi\Vert^2$ and treat $F$ as a function 
of the complex coefficients $a^1,\,\ldots,a^m$ from \mythetag{2.1}, then 
the minimum of the function $F(a^1,\,\ldots,a^m)$ is determined by the
vanishing conditions for its partial derivatives:
$$
\hskip -2em
\frac{\partial F}{\partial a^i}=0,\quad
\frac{\partial F}{\partial\overline{a^i}}=0
\text{, \ where \ }i=1,\,\ldots,\,n. 
\mytag{2.6}
$$
Calculating the derivatives \mythetag{2.6} for \mythetag{2.1}, we derive
the equations 
$$
\xalignat 2
&\hskip -2em
\sum^n_{i=1}g_{ij}\,\overline{a^i}=\bigl<f\bigl|\kern 1pt e^{\kern 0.5pt\lambda_jx}\bigr>,
&&\sum^n_{j=1}g_{ij}\,a^j=\bigl<e^{\kern 0.5pt\lambda_ix}\bigl|f\bigr>,
\mytag{2.7}
\endxalignat
$$
where $i=1,\,\ldots,\,n$. Due to \mythetag{2.5} two sets of equations 
\mythetag{2.7} differ from each other only by complex conjugation. 
Hence it is sufficient to solve only one of them.\par 
     We choose for solving the second set of the equations \mythetag{2.7}.
It is solved with the use of the inverse Gram matrix $G^{-1}$. Let's
denote through $g^{ij}$ the components of the transpose of the inverse
Gram matrix $(G^{-1})^{\sssize\top}$. Then the quantities $g_{ij}$ and
$g^{ij}$ are related to each other as follows:
$$
\xalignat 2
&\hskip -2em
\sum^n_{k=1}g_{ik}\,g^{jk}=\delta^j_i,
&&\sum^n_{k=1}g^{kj}\,g_{ki}=\delta^j_i.
\mytag{2.8}
\endxalignat
$$ 
Here $\delta^j_i$ are the components of the unit matrix. They are called 
Kronecker's delta.\par
     Applying the second relationship \mythetag{2.8} to the second set of
equations \mythetag{2.7}, we get their solution. This solution is given by 
the formula 
$$
\hskip -2em
a^j=\sum^n_{i=1} g^{ij}\,\bigl<e^{\kern 0.5pt\lambda_ix}\bigl|f\bigr>\text{, \ where \ }
j=1,\,\ldots,\,n. 
\mytag{2.9}
$$
Substituting \mythetag{2.9} into \mythetag{2.1}, we derive 
$$
\hskip -2em
F_{\text{min}}=\Vert f\Vert^2-\sum^n_{i=1}\sum^n_{j=1}g^{ij}\,
\bigl<e^{\kern 0.5pt\lambda_ix}\bigl|f\bigr>\,\bigl<f\bigl|\kern 0.5pt 
e^{\kern 0.5pt\lambda_jx}\bigr>.
\mytag{2.10}
$$
The formulas  \mythetag{2.9} and \mythetag{2.10} yield a solution of the
linear approximation problem.\par
\head
3. One frequency approximation for the sign function. 
\endhead
     The minimum value $F_{\text{min}}$ in \mythetag{2.10} is a function of
$\lambda_1,\,\ldots,\,\lambda_n$. Let's denote
$$
\hskip -2em
\Phi=\Phi(\lambda_1,\,\ldots,\lambda_n)=F_{\text{min}}.
\mytag{3.1}
$$
The nonlinear approximation problem consists in finding the absolute minimum
of the function \mythetag{3.1} as $\lambda_1,\,\ldots,\,\lambda_n$ run over
$\Bbb C^n$. Though potentially this could be not a minimum, but infimum. In 
any case, since $0\leqslant F_{\text{min}}\leqslant\Vert f\Vert^2$, the minimal 
value of the function $\Phi(\lambda_1,\,\ldots,\lambda_n)$ does exist and is 
finite.\par
     Let's consider the case $n=1$. We call it the one frequency case since the
quantities $\lambda_1,\,\ldots,\,\lambda_n$ are often associated with 
eigenfrequencies in applications. In order to study this case thoroughly we
choose 
$$
\hskip -2em
f(x)=\sign(x)=\cases -1&\text{\ \ for \ }x<0,\\
1&\text{\ \ for \ }x\geqslant 0
\endcases
\mytag{3.2}
$$
as an example. In the one frequency case the formula \mythetag{2.10} simplifies. 
It turns to 
$$
\hskip -2em
F_{\text{min}}=\Vert f\Vert^2-\frac{\bigl<e^{\kern 0.5pt\lambda_1x}\bigl|f\bigr>
\,\bigl<f\bigl|\kern 0.5pt e^{\kern 0.5pt\lambda_1x}\bigr>}{g_{11}}.
\mytag{3.3}
$$\par
     The norm of the function \mythetag{3.2} is easily calculated: $\Vert f\Vert=1$.
The denominator $g_{11}$ in the formula \mythetag{3.3} is also easily calculated:
$$
\hskip -2em
g_{11}=\frac{1}{2\,\pi}\!\int\limits^{\,+\pi}_{\!-\pi}\!\!
e^{\kern 0.5pt(\lambda_1+\bar\lambda_1)x}\,d\kern 0.5pt x=
\frac{e^{\kern 0.5pt(\lambda_1+\bar\lambda_1)\pi}-
e^{-(\lambda_1+\bar\lambda_1)\pi}}{2\,\pi\,(\lambda_1+\bar\lambda_1)
\vphantom{\vrule height 10pt depth 0pt}}.
\mytag{3.4}
$$
Now let's calculate the quantities $\bigl<e^{\kern 0.5pt\lambda_1x}\bigl|f\bigr>$
and $\bigl<f\bigl|\kern 0.5pt e^{\kern 0.5pt\lambda_1x}\bigr>$ in \mythetag{3.3}:
$$
\align
&\hskip -2em
\bigl<e^{\kern 0.5pt\lambda_1x}\bigl|f\bigr>=
\frac{1}{2\,\pi}\!\int\limits^{\,+\pi}_{\!-\pi}\!\!
e^{\kern 0.5pt\bar\lambda_1x}\sign(x)\,d\kern 0.5pt x
=\frac{e^{\kern 0.5pt\bar\lambda_1\pi}+e^{-\bar\lambda_1\pi}}{2\,\pi\,\bar\lambda_1
\vphantom{\vrule height 10pt depth 0pt}}-\frac{1}{\pi\,\bar\lambda_1
\vphantom{\vrule height 10pt depth 0pt}},
\mytag{3.5}\\
&\hskip -2em
\bigl<f\bigl|\kern 0.5pt e^{\kern 0.5pt\lambda_1x}\bigr>=
\frac{1}{2\,\pi}\!\int\limits^{\,+\pi}_{\!-\pi}\!\!
\sign(x)\,e^{\kern 0.5pt\lambda_1x}\,d\kern 0.5pt x
=\frac{e^{\kern 0.5pt\lambda_1\pi}+e^{-\lambda_1\pi}}{2\,\pi\,\lambda_1
\vphantom{\vrule height 10pt depth 0pt}}-\frac{1}{\pi\,\lambda_1
\vphantom{\vrule height 10pt depth 0pt}}.
\mytag{3.6}
\endalign
$$
The spectral parameter $\lambda_1$ is a complex variable. Therefore we
write
$$
\hskip -2em
\lambda_1=u+\frak i\,v\text{, \ where \ }\frak i=\sqrt{-1\kern 0.5pt}.
\mytag{3.7}
$$
Substituting \mythetag{3.7} into the formula \mythetag{3.4}, we derive 
$$
g_{11}=\frac{\sinh(2\,\pi\,u)}{2\,\pi\,u}.
$$
Similarly, using the formulas \mythetag{3.5} and \mythetag{3.6}, we obtain
$$
\gathered
\bigl<e^{\kern 0.5pt\lambda_1x}\bigl|f\bigr>
\bigl<f\bigl|\kern 0.5pt e^{\kern 0.5pt\lambda_1x}\bigr>
=\frac{\cosh(2\,\pi\,u)+\cos(2\,\pi\,v)}{2\,\pi^2(u^2+v^2)}+\\
+\frac{2-4\,\cosh(\pi\,u)\,\cos(\pi\,v)}{2\,\pi^2(u^2+v^2)}
=\frac{(\cosh(\pi\,u)-\cos(\pi\,v))^2}{\pi^2(u^2+v^2)}.
\endgathered
$$
As a result we obtain the following expression for the function \mythetag{3.1}:
$$
\hskip -2em
\Phi(\lambda_1)=\Phi(u+\frak i\,v)=1-\frac{2\,u}{u^2+v^2}\,
\frac{(\cosh(\pi\,u)-\cos(\pi\,v))^2}
{\pi\,\sinh(2\,\pi\,u)}. 
\mytag{3.8}
$$
Passing to the limit $u\to 0$ in \mythetag{3.8}, we obtain the function
$$
\hskip -2em
\Phi(0+\frak i\,v)=1-\frac{(1-\cos(\pi\,v))^2}{\pi^2\,v^2}. 
\mytag{3.9}
$$
The function \mythetag{3.9} has two absolute minima at $v=\pm\kern 0.5pt v_0$, where 
$v_0$ is a real irrational number being a solution of the equation
$$
\hskip -2em
\sin(\pi\,v)\,\pi\,v+\cos(\pi\,v)=1
\mytag{3.10}
$$
and such such that $0.1<v_0<0.9$. Solving \mythetag{3.10} numerically, we find that
$$
\pagebreak 
\hskip -2em
v_0=0.742019\ldots .
\mytag{3.11}
$$
Applying \mythetag{3.10} and \mythetag{3.11} to \mythetag{3.9}, we derive
$$
\hskip -2em
\Phi_{\text{min}}=\Phi(\pm\kern 0.5pt \frak i\,v_0)=\cos^2(v_0)=
0.4749383\ldots .
\mytag{3.12}
$$
One can show that two complex numbers $\lambda=\pm\kern 0.5pt\frak i\,v_0$ are 
absolute minima of the function \mythetag{3.8} as well. The number \mythetag{3.12} 
is its minimum value.
\par
    It is curious to note that the best one frequency root mean square approximation
for the real-valued function $\sign(x)$ is given by a complex function, which is 
quite unlike to its Fourier expansion approximation. 
\head
4. Two frequencies approximation for the sign function. 
\endhead
     The two frequencies case $n=2$ subdivides into two subcases $\lambda_1\neq
\lambda_2$ and $\lambda_1=\lambda_2$. The subcase $\lambda_1=\lambda_2$ is 
a cluster case (see \mycite{2}). In this case two exponential functions 
$e^{\kern 0.5pt\lambda_1x}$ and $e^{\kern 0.5pt\lambda_2x}$ are replaced by  
two expo-polynomials
$$
\xalignat 2
&\phi_1(x)=e^{\kern 0.5pt\lambda_1x},
&&\phi_2(x)=x\,e^{\kern 0.5pt\lambda_1x}.
\endxalignat
$$
The formula \mythetag{2.3} is replaced by the formula
$$
g_{ij}=\bigl<\phi_i\bigl|\phi_j\bigr>.
$$
Similarly, the formula \mythetag{2.10} is replaced by the formula
$$
\hskip -2em
F_{\text{min}}=\Vert f\Vert^2-\sum^n_{i=1}\sum^n_{j=1}g^{ij}\,
\bigl<\phi_i\bigl|f\bigr>\,\bigl<f\bigl|\phi_j\bigr>.
\mytag{4.1}
$$
The formula \mythetag{3.1} remains unchanged.
\par
     For the beginning we consider the cluster case $\lambda_1=\lambda_2$ with 
$\lambda_1=0+\frak i\,v$. Applying the formulas \mythetag{3.1} and \mythetag{4.1} 
to this case, we derive  
$$
\hskip -2em
\gathered
\Phi(0+\frak i\,v)=1-(1-\cos(\pi\,v))\,\times\\
\vspace{1ex}
\times\,\frac{(2\,\cos(\pi\,v)\,\pi^2\,v^2-3\,\cos(\pi\,v)+3
+4\,\pi^2\,v^2-6\,\sin(\pi\,v)\,\pi\,v)}{\pi^4\,v^4}.
\endgathered
\mytag{4.2}
$$
Unlike \mythetag{3.9}, the function \mythetag{4.2} has exactly one absolute 
minimum at $v=0$ (i\.\,e\. at the origin of the complex plane $\lambda_1=0
+\frak i\,0=0$) such that
$$
\hskip -2em
\Phi_{\text{min}}=\lim_{v\to\,0}\Phi(0+\frak i\,v)=\frac{1}{4}.
\mytag{4.3}
$$
.\par
     The general cluster case corresponds to $\lambda_1=\lambda_2$ with 
$\lambda_1=u+\frak i\,v$. In this case, applying \mythetag{3.1} 
and \mythetag{4.1}, we obtain an explicit expression for the function
$\Phi(\lambda_1)=\Phi(u+\frak i\,v)$. However, this expression is much more
bulky than the expression \mythetag{3.8}. Analyzing this bulky expression
numerically, we find that the function $\Phi(\lambda_1)$ has a unique
minimum at the same point $\lambda_1=0$ as the function \mythetag{4.2}.
\par
     The next step is to proceed to the non-cluster case. In this case the
formulas $\mythetag{2.10}$ and $\mythetag{3.1}$ yield a function of two
complex variables $\Phi(\lambda_1,\lambda_2)$. The expression for 
$\Phi(\lambda_1,\lambda_2)$ is very bulky and complicated. \pagebreak 
Finding an absolute minimum numerically for such a function is also a 
complicated problem. Therefore at present moment we can only formulate 
a conjecture that
$$
\hskip -2em
\lim_{\Sb\lambda_1\to \,0\\\lambda_2\to\,0\endSb}\Phi(\lambda_1,\lambda_2)=
\inf_{\lambda_1\neq\lambda_2}\{\Phi(\lambda_1,\lambda_2)\}=\frac{1}{4}.
\mytag{4.4}
$$
The formula \mythetag{4.3} provides a support for our conjecture \mythetag{4.4}. 
\head
5. Conclusions. 
\endhead
     Once a square integrable function $f(x)\in L^2([\kern 1pt-\pi,\,+\pi])$
is given, the optimal values of $\lambda_1,\,\ldots,\,\lambda_n$ for $f$
are called an $n$-frequencies spectrum of the function $f$. The above example
of the sign function shows that the $n$-frequencies spectrum is not unique.
Moreover, it is not stable. The spectral point $\lambda_1=\frak i\,v_0$,
which is present in one frequency spectrum, disappears in two frequencies 
spectrum.\par

\enddocument
\end